\documentclass[12pt]{article}
\usepackage{float} 
\usepackage{amsmath}
\usepackage{amssymb}
\usepackage{a4}
\usepackage{theorem}
\usepackage{epsfig}
\usepackage[mathscr]{eucal}

\newcommand{\bzero}{\mathbf{0}}
\newcommand{\mlt}{\prec}
\newcommand{\mgt}{\succ}
\newcommand{\mlet}{\preccurlyeq}
\newcommand{\mget}{\succcurlyeq}
\newcommand{\mlth}{\prec_h}
\newcommand{\mgth}{\succ_h}

\newcommand{\fele}{f_1,\ldots,f_r}
\newcommand{\fset}{{\{\fele}\}}
\newcommand{\balpha}{{\pmb{\alpha}}}
\newcommand{\bbeta}{{\pmb{\beta}}}
\newcommand{\bgamma}{{\pmb{\gamma}}}
\newcommand{\bdelta}{{\pmb{\delta}}}
\newcommand{\btau}{{\pmb{\tau}}}
\newcommand{\txi}{\tilde{\xi}}

\newcommand{\tA}{\tilde{A}}
\newcommand{\CC}{\mathbb{C}}

\newcommand{\sF}{\mathscr{F}}

\newcommand{\tf}{\tilde{f}}
\newcommand{\oI}{\overline{I}}

\newcommand{\tj}{\tilde{j}}
\newcommand{\tk}{\tilde{k}}
\newcommand{\sL}{\mathscr{L}}
\newcommand{\sM}{\mathscr{M}}

\newcommand{\NN}{\mathbb{N}}

\newcommand{\cO}{{\cal O}}

\newcommand{\QQ}{\mathbb{Q}}
\newcommand{\RR}{\mathbb{R}}

\newcommand{\sT}{\mathscr{T}}

\newcommand{\bx}{\mathbf{x}}
\newcommand{\bxa}{\bx^\balpha}
\newcommand{\bxb}{\bx^\bbeta}
\newcommand{\bxtb}{\bx^{\tilde{\bbeta}}}
\newcommand{\bxg}{\bx^\bgamma}
\newcommand{\bxt}{\bx^\btau}

\newcommand{\Kx}{K[\bx]}
\newcommand{\Ktx}{K[t,\bx]}
\newcommand{\loc}{{\mlt}}
\newcommand{\Kxloc}{S_\loc^{-1}\Kx}
\newcommand{\Iloc}{S_\loc^{-1}I}

\newcommand{\NormalForm}{{\rm\bf NormalForm}}
\newcommand{\StandardBasis}{{\rm\bf StandardBasis}}
\newcommand{\StandardBasisChange}{{\rm\bf StandardBasisChange}}
\newcommand{\ReducingSet}{{\rm\bf ReducingSet}}
\newcommand{\GaussEliminate}{{\rm\bf GaussEliminate}}
\newcommand{\ComputeZgen}{{\rm\bf ComputeZgen}}

\newcommand{\lcm}{\operatorname{lcm}}

\newcommand{\spoly}{\operatorname{spoly}}
\newcommand{\Spoly}{\operatorname{Spoly}}

\newcommand{\sa}{\hspace*{1cm}}
\newcommand{\saa}{\sa\sa}
\newcommand{\saaa}{\sa\sa\sa}

{\theorembodyfont{\rm}
    \newtheorem{definition}{Definition}[section]
    
    \newtheorem{example}[definition]{Example}
    \newtheorem{algorithmxx}[definition]{Algorithm}
}
\newtheorem{lemma}[definition]{Lemma}
\newtheorem{corollary}[definition]{Corollary}

\newtheorem{proposition}[definition]{Proposition}
\newtheorem{theorem}[definition]{Theorem}
\newcommand{\proof}{{\noindent\it Proof:\ }}

\newcommand{\proofend}{$\square$\medskip\par}
\newenvironment{algorithm}{\samepage
\begin{algorithmxx}\ \newline
\rule{\linewidth}{0.3mm}
\begin{minipage}{\linewidth}}{\end{minipage}\newline
\rule{\linewidth}{0.3mm}
\end{algorithmxx}}
\newcommand{\bibauthor}[2]{{\it{#2}~{#1},}}
\newcommand{\bibtitlea}[1]{#1,}            
\newcommand{\bibtitleb}[1]{{\it #1},}      
\newcommand{\bibcompany}[1]{#1}

\newcommand{\bibyear}[1]{(#1),}
\newcommand{\bibyearx}[1]{(#1)}
\newcommand{\bibjournal}[2]{#1\ {\bf #2}}
\newcommand{\bibpages}[2]{{#1}--{#2}}
\newcommand{\bibend}{.}

\begin{document}
\newcounter{hlp}
\ \vspace*{5cm}\\
\noindent{\LARGE\bf Standard bases with respect to
the Newton filtration\\[2ex]}
{\large\bf Stephan Endra\ss\\[1ex]}
Johannes Gutenberg-Universit\"at, Mainz, Germany\\[2ex]
\today\\[5ex]
{\bf Abstract}\\[1ex]
{\small The aim of this article is to introduce standard bases
of ideals in polynomial rings
with respect to a class of orderings which are not
necessarily semigroup orderings. Our approach generalises
the concept of standard bases with respect to semigroup
orderings described in \cite{graebe,greuelpfister}.
To compute these standard bases we give a slightly
modified version of the Buchberger algorithm \cite{buchberger}.
The orderings we consider are refinements of certain filtrations.
In the local case these filtrations are Newton filtrations.
For a zero dimensional ideal, an algorithm converting
standard bases with respect to local orderings is given.
As an application, we show how to compute the spectrum of an
isolated complex hypersurface singularity
$f\colon(\CC^n,\bzero)\rightarrow(\CC,0)$
with nondegenerate principal part.}
\setcounter{section}{-1}
\section{Introduction}
The theory of standard bases has developed since B.~Buchberger
introduced standard bases of ideals in polynomial rings
with respect to semigroup wellorderings
(now called Gr\"obner bases) in 1965 \cite{buchberger}.
It has been extended to certain localisations of
the polynomial rings by allowing semigroup orderings
which are not wellorderings \cite{mora,graebe,greuelpfister,moraseven}.
The main tool to compute such a standard basis is the
Buchberger algorithm \cite{buchberger,beckerweispfenning}.
By choosing an appropriate
semigroup ordering, invariants of ideals may be computed
from a standard basis with respect to this ordering.
For example, the Hilbert function can be obtained from
a standard basis with respect to a degree ordering.
However there are some invariants which are not
related to semigroup orderings. Consider an isolated
hypersurface singularity $f\colon\CC^n\rightarrow\CC$
with nondegenerate principal part. Its Milnor number
$\mu(f)$ is readily computed from a standard basis of
the Jacobian ideal $J_f$ of $f$ with respect to any
local semigroup ordering (in fact with respect
to any local ordering as will be shown). However 
the spectrum of $f$ is the Poincare series of
$\cO_{\CC^n,\bzero}/J_f$ graded with respect to the
Newton filtration given by $f$ \cite{saito,kv}.
The Newton filtration can be refined to a local
semigroup ordering if and only if $f$ is semiquasihomogeneous.
Here we study standard bases with respect to orderings
which are refinements of Newton filtrations. In general, such
an ordering is not a semigroup ordering.

In the first section normal and noetherian orderings are introduced.
This is the class of orderings for which we are able to
give a normal form algorithm which terminates
(algorithm \ref{algorithm:normalform}).
The setup for
standard bases with respect to normal orderings is outlined
in the second section.
To compute standard bases, $s$-polynomial sets
and reducing sets are introduced. With their help a
modified Buchberger algorithm can be formulated 
(algorithm \ref{algorithm:standardbasis}).
It terminates for any normal noetherian ordering for
which reducing sets exist (theorem \ref{theorem:standardbasis})
and returns a standard basis (theorem \ref{theorem:standardbasisx}).

In the next section Newton orderings are introduced as
refinements of Newton filtrations.
Some Newton orderings are not normal (lemma \ref{lemma:newtonnormal}),
but they all are
noetherian and admit reducing sets
(propositions \ref{proposition:newtonnoether}
and \ref{proposition:newtonred}).
The case of a zero dimensional ideal and a local ordering is
studied in in the fourth section. Here standard bases
do always exist (proposition \ref{proposition:always}).
Moreover a standard basis with respect to a given local ordering
can be computed from a standard basis with respect to another
local ordering using only linear algebra
(algorithm \ref{algorithm:convert}).
In the last section we show how to compute the
spectrum of an isolated complex hypersurface singularity
with nondegenerate principal part (corollary \ref{corollary:spectrum}).
The author implemented the computation of the spectrum into the computer
algebra program {\tt Singular} \cite{singular} using algorithm \ref{algorithm:convert}.
\section{Orderings}\label{sect:orderings}
First fix some notation.
Let $\NN=\{0,1,2,\ldots\}$ be the nonnegative integers,
$K$ a field and $\Kx=K[x_1,\ldots,x_n]$ the polynomial ring in
$n\geq 1$ indeterminates over $K$.
We use the exponent notation $\bxa=
x_1^{\alpha_1}\cdots x_n^{\alpha_n}$,
$\balpha=(\alpha_1,\ldots,\alpha_n)\in\NN^n$.
Let
\begin{align*}
    \sM &= \{\bxa\mid\balpha\in\NN^n\}\text{\ be the set of all
           monomials of $\Kx$ and}\\
    \sT &= \{c\,\bxa|c\in K^\ast, \bxa\in\sM\}\text{\ the set of all
           terms in $\Kx$.}
\end{align*}
The polynomials of $\Kx$ are sums of the form
$\sum_{\balpha\in A}c_\balpha\bxa$, where $A\subset\NN^n$ is
a finite set and $c_\balpha\bxa\in\sT$ for all $\balpha\in A$
(we do not write monomials with zero coefficient).
Consider a total ordering $\mlt$ on $\NN^n$ (which is not not necessarily
a semigroup ordering) and denote the induced ordering on $\sM$ also by $\mlt$.
\begin{definition}
    Let $G=\fset\subset\Kx$ be a finite set of polynomials,
    $I\subseteq\Kx$ an ideal and
    $f=\sum_{\balpha\in A}c_\balpha\bxa\in \Kx$.
    \begin{itemize}
        \item[1)] $LM(f)=\bx^{\max A}$ is called the {\em lead monomial of $f$}.
        \item[2)] $LC(f)=c_{\max A}$ is called the
                  {\em lead coefficient of $f$}.
        \item[3)] $LT(f)=LC(f)\cdot LM(f)$ is called the
                  {\em lead term of $f$}.
        \item[4)] $L(G)=\{LM(\bxa g)\mid \bxa\in\sM,
                  g\in G\}$ is called the
                  {\em lead monomial set of $G$}.
        \item[5)] $L(I)=\{LM(g)\mid g\in I\}$ is called the
                  {\em lead monomial set of $I$}.
    \end{itemize}
\end{definition}
Any nonempty set of monomials which is the lead monomial
set of a finite set of polynomials or of an ideal will be called
{\em lead monomial set}.
For a semigroup ordering $L(I)$ and
$L(G)$ can be identified with monomial ideals.
In general this is not the case.
But the ordering $\mlt$ still satisfies
\begin{equation*}
    LM(f+g)\mlet\max\left( LM(f),LM(g)\right)
\end{equation*}
for all $f$, $g\in \Kx$. Instead of respecting the semigroup
structure of $\NN^n$, the ordering $\mlt$ should satisfy the
following condition.
\begin{definition}
    The ordering $\mlt$ is called {\em normal}, if
    for all $f\in \Kx$,
    $\bxa\in\sM$
    \begin{align*}
        %
        \bxa\mlt 1 \Longrightarrow LM(\bxa f) \mlt LM(f).
    \end{align*}
    A normal ordering $\mlt$ is called {\em global} (resp.~{\em local})
    if $x_i\mgt 1$ (resp.~$x_i\mlt 1$), $1\leq i\leq n$. A normal
    ordering is called {\em mixed} if it is neither global
    nor local.
\end{definition}
For a normal ordering the set
$S_\loc=\{g\in\Kx\mid LM(g)=1\}$
is multiplicatively closed. Moreover $LM(f)=LM(gf)$ for all
$f$, $g\in\Kx$ with $LM(g)=1$.
For semigroup orderings the lead monomial sets
behave like ideals (in a noetherian ring), i.e.~every
increasing sequence gets stationary.
\begin{definition}
    The ordering $\mlt$ is called {\em noetherian} if the following
    condition is satisfied:
    \begin{itemize}
        \item[] Every increasing sequence of lead monomial sets
                $L_0\subseteq L_1\subseteq L_2\subseteq\ldots$
                gets stationary.
    \end{itemize}
\end{definition}
Any semigroup ordering is normal and noetherian, but the converse
does not hold.
\section{Standard bases}\label{sect:general}
Proceeding along the lines of \cite{greuelpfister}, we introduce
standard bases with respect to normal orderings.
From now on let $\mlt$ be a normal ordering.
Let $\Kxloc$ be the localisation of $\Kx$ in $S_\loc$. For every ideal
$I\subseteq\Kx$, let $\Iloc = I\otimes_{\Kx}\Kxloc$ and
$\oI=\Iloc\cap\Kx$. Note that always $L(I)=L(\oI)$.
\begin{definition}\label{definition:standard}
    Let $I\subseteq\Kx$ be an ideal and $G=\fset\subseteq\oI$.
    \begin{itemize}
        \item[1)] $G$ is called a
            {\em standard basis of $I$} if $L(G)=L(I)$.
        \item[2)] $G$ is called {\em interreduced}
            if $L(\{f_i\})\not\subseteq L(G\setminus\{f_i\})$
                for all $i=1,\ldots,r$.
        \item[3)] $G$ is called {\em reduced}
            if for every $f\in G$ no monomial of $f$ except
            its lead monomial is contained in $L(G)$.
    \end{itemize}
\end{definition}
If $\mlt$ is not a semigroup ordering, then in general not every
ideal has a standard basis. Albeit for noetherian orderings this
is true.
\begin{lemma}
    If $\mlt$ is noetherian, the every ideal $I\subseteq\Kx$
    has a standard basis.
\end{lemma}
\proof The standard proof applies.
\proofend
Every standard basis $G$ of $I$ can be shortened to an interreduced
standard basis by iteratively deleting those $f\in G$
for which $L(\{f\})\subseteq L(G\setminus\{f\})$. Reduced standard bases do
in general only exist for global orderings.
To do standard basis computations a normal form is needed:
\begin{definition}\label{definition:normal_form}
    \cite[def.~1.5]{greuelpfister} Let $\sF=\{G\subset \Kx\mid
    G\text{\ ordered and finite}\}$. A function
    $NF\colon \Kx\times\sF\rightarrow \Kx$ is called
    {\em normal form} if for all $p\in \Kx$ and
    $G\in\sF$ $NF(p,G)\neq 0 \Rightarrow LM(NF(p,G))\not
    \in L(G)$. Then $NF(p,G)$ is called {\em normal form
    of $p$ with respect to $G$}.
\end{definition}
Following \cite{greuelpfister}, to get a normal form which
does not only work in the global case, one considers homogeneous
polynomials in $K[t,\bx]$. So let
\begin{align*}
    \sM_t &=\{t^{\alpha_0}\bxa\mid \alpha_0\in\NN,
            \balpha\in\NN^n\},\\
    \sT_t &=\{c\,t^{\alpha_0}\bxa\mid c\in K^\ast,
            \alpha_0\in\NN, \balpha\in\NN^n\},
\end{align*}
be the monomials and terms of $\Ktx$.
Define a global wellordering $\mlth$ on $\sM_t$ by setting
\begin{align*}
    t^{\alpha_0}\bxa\mlth t^{\beta_0}\bxb\Longleftrightarrow &
    \ \alpha_0 +|\balpha| < \beta_0+|\bbeta|\quad\text{or}\\
    &\ \alpha_0 +|\balpha| = \beta_0+|\bbeta|\quad\text{and}\quad\bxa\mlt\bxb.
\end{align*}
With respect to powers of $t$, this ordering behaves like a
semigroup ordering.
We frequently need to consider lead monomials in $\Ktx$ and $\Kx$.
For every polynomial $f=\sum_{\balpha\in A}c_\balpha\bxa\in \Kx$,
let $f^h=\sum_{\balpha\in A}c_\balpha t^{\deg f-|\balpha|}\bxa$
be the homogenisation with respect to $t$.
Conversely, for $f=\sum_{(\alpha_0,\balpha)\in A}
c_{(\alpha_0,\balpha)}t^{\alpha_0}\bxa$ let
$\left.f\right|_{t=1}=
\sum_{(\alpha_0,\balpha)\in A}
c_{(\alpha_0,\balpha)}\bxa$ be its dehomogenisation.
For polynomials of $\Kx$ (resp.~$\Ktx$) we always use
$\mlt$ (resp.~$\mlth$) in the computation of lead monomials,
lead terms, \ldots .
Let $G=\fset\subset\Ktx$ be a finite ordered set of
homogeneous polynomials and let $p\in\Ktx$ homogeneous.
For a semigroup ordering, the following is just the normal form
{\bf NFMora} of \cite{greuelpfister}.
%
%
\begin{algorithm}\label{algorithm:normalform}
    $h:=\NormalForm(p,G)$\\
    \sa $h:=p$\\
    \sa $H:=\emptyset$\\
    \sa WHILE $\left(\begin{array}{l}
              \text{exist $f\in G\cup H$, $\eta\in\sT_t$,
              $\alpha\in\NN$ with}\\
              \text{$LT(\eta f)=LT(t^\alpha h)$
              and $\eta\mlet 1$ if $f\in H$}
              \end{array}\right)$ DO\\
    \saa choose first such $f$ with $\alpha$ minimal\\
    \saa IF $\alpha>0$ THEN\\
    \saaa $H:=H\cup\{h\}$\\
    \saa $h:=t^\alpha h-\eta f$\\
    \saa IF $t\mid h$ THEN\\ 
    \saaa choose $\alpha$ maximal with $t^\alpha\mid h$\\
    \saaa $h:=h/t^\alpha$
\end{algorithm}
\begin{lemma}\label{lemma:normal_form}
    If $h$ is a normal form of $p$ with respect
    to $G=\fset$ computed by \NormalForm, then there
    exist $\xi_1,\ldots,\xi_k\in\sT_t$, $g\in\Ktx$ homogeneous and
    $j_1,\ldots,j_k\in\{1,\ldots,r\}$
    such that
    \begin{equation*}
        gp=\sum_{i=1}^k\xi_if_{j_i} + h
    \end{equation*}
    with $LM(g|_{t=1})=1$ and $LM(h|_{t=1})\not\in L(G|_{t=1})$.
    Moreover $gp$, $\xi_1 f_{j_1},\ldots,\xi_k f_{j_k}$ and
    $h$ (if nonzero) are homogeneous of the same degree with
    \begin{equation*}
        LM(gp) = LM(\xi_1f_{j_1})\mgth\ldots\mgth
        LM(\xi_k f_{j_k})\mgth LM(h).
    \end{equation*}
\end{lemma}
\proof This just the proof of \cite[theorem 1.9 2)]{greuelpfister}.
\proofend
\begin{proposition}\label{proposition:normal_form}
    If the ordering $\mlt$ is noetherian, then
    \NormalForm\ is a normal form in the sense
    of definition \ref{definition:normal_form}.
\end{proposition}
\proof In view of lemma \ref{lemma:normal_form} it is enough to show that
\NormalForm\ terminates. But $\mlt$ is noetherian, so 
the proof of \cite[prop.~1.9 1)]{greuelpfister} applies.
\proofend
To actually compute such a standard basis one needs $s$-polynomials.
\begin{definition}
    Let $f=\sum_{\balpha\in A}c_\balpha\bxa$,
    $g=\sum_{\bbeta\in B}d_\bbeta\bxb\in \Kx$.
    Then
    \begin{equation*}
        \spoly (f,g)_{(\balpha,\bbeta)}=
        \left(d_\bbeta\bxb f-c_\balpha\bxa g\right)/\gcd\left(\bxa,\bxb\right)
    \end{equation*}
    is called {\em $s$-polynomial of $f$ and $g$ at $(\balpha,\bbeta)$}.
    Moreover
    \begin{align*}
        \Spoly(f,g)=\left\{
        \spoly(f,g)_{(\balpha,\bbeta)}\mid
        \begin{array}{l}
        (\balpha,\bbeta)\in A\times B, \exists \bx^\bgamma\in\sM:\\
        LM(\bx^\bgamma\spoly(f,g)_{(\balpha,\bbeta)})=
        \bx^\bgamma\lcm(\bxa,\bxb)
        \end{array}
        \right\}
    \end{align*}
    is called the {\em $s$-polynomial set of $f$ and $g$}.
\end{definition}
For a semigroup ordering 
$\left|\Spoly(f,g)\right|=1$ and we get the usual
$s$-polynomial.
Now we want to imitate the Buchberger algorithm.
Given a finite set of generators
$G=\fset$ of an ideal $I$, this algorithm
enlarges $G$ by some elements of $I$ such that all $s$-polynomials
$\spoly(f,g)$, $f,g\in G$ reduce to zero. For a semigroup ordering
this implies that $\bxa\spoly(f,g)$ reduces to zero
for all $\bxa\in\sM$. In our setup this is not the case.
Let $G=\fset\subset\Ktx$ be a finite set of homogeneous
polynomials and let $f\in I=(\fele)$ homogeneous.
Assume that that \NormalForm\ terminates for $\mlt$.
\begin{definition}\label{definition:reducing_set}
    A finite set $R(f,G)\subset\sM$ of monomials is called a
    {\em reducing set for $(f,G)$} if the following holds:
    For all $\bxa\in\sM$, there exist 
    \begin{itemize}
        \item a homogeneous polynomial $g\in\Ktx$ with $LM(g|_{t=1})=1$,
        \item $\xi_1,\ldots,\xi_k\in\sT_t$ and
        \item $h_1,\ldots,h_k\in G\cup\{\NormalForm(\bxb f,G)\mid
              \bxb\in R(f,G)\}$
    \end{itemize}
    such that the polynomial $\bxa f$ has
    a representation of the form
    \begin{equation*}
        g\bxa f = \sum_{i=1}^k\xi_i h_i
    \end{equation*}
    with $LM(g\bxa f)=LM(\xi_1 h_1)\mgth LM(\xi_i h_i)$, $2\leq i\leq k$.
    If there exists a reducing set for every $(f,G)$ as
    above, then we say that {\em reducing sets exist}.
\end{definition}
For a semigroup ordering we always can take $R(f,G)=\{1\}$.
Now we can formulate a standard basis algorithm. Let
$G=\fset\subset\Kx$ be a finite set of polynomials.
Let $\mlt$ be such that \NormalForm\ terminates and that
reducing sets exist. Moreover assume that the procedure
$\ReducingSet(h,S)$ computes a reducing set for $(h,S)$.
\par
\begin{algorithm}\label{algorithm:standardbasis}\ \newline\noindent
    $S:=\StandardBasis(G)$\\
    \hspace*{1cm}$S:=G^h$\\
    \hspace*{1cm}$P:=\{(f,g)\mid f,g\in S\}$\\
    \hspace*{1cm}WHILE $P\neq\emptyset$ DO\\
    \hspace*{2cm}choose $(f,g)\in P$; $P:=P\setminus\{(f,g)\}$\\
    \hspace*{2cm}FOR ALL $h\in\Spoly(f,g)$ DO\\
    \hspace*{3cm}R:=\ReducingSet(h,S)\\
    \hspace*{3cm}FOR ALL $\bxa\in R$ DO\\
    \hspace*{4cm}$p:=\NormalForm(\bxa h,S)$\\
    \hspace*{4cm}IF $p\neq 0$ THEN\\
    \hspace*{5cm}$S:=S\cup\{p\}$\\
    \hspace*{5cm}$P:=P\cup\{(p,f)\mid f\in S\}$\\
    \hspace*{1cm}$S:=S|_{t=1}$
\end{algorithm}
\begin{theorem}\label{theorem:standardbasis}
    Let $\mlt$ be noetherian such that
    reducing sets exist. Then \StandardBasis\ terminates.
\end{theorem}
\proof Again the standard proof applies.
\proofend
\begin{theorem}\label{theorem:standardbasisx}
    Let $\mlt$ be such that
    \NormalForm\ terminates and reducing sets exist. Let
    $I\subseteq\Kx$ be an ideal.
    Equivalent for $G=\fset\subseteq\oI$ are:
    \begin{itemize}
        \item[  i)] $G$ is a standard basis of $I$.
        \item[ ii)] $G=\StandardBasis(G)$.
        \item[iii)] $\NormalForm(\bxa g^h,G^h)=0$ for all
                    $\bxa\in\sM$, $h\in\Spoly(f_i,f_j)$,
                    $1\leq i\leq j\leq r$.
        \item[ vi)] $f\in\oI\Leftrightarrow\NormalForm(f^h,G^h)=0$.
        \item[  v)] $f\in\oI\Leftrightarrow gf=\sum_{i=1}^k\xi_i f_{j_i}$
                    for some
                    $g\in\Kx$,
                    $\xi_1,\ldots,\xi_k\in\sT$ and
                    $j_1,\ldots,j_k\in\{1,\ldots,r\}$ with
                    $LM(g)=1$ and 
                    $LM(f) = LM(\xi_1 f_{j_1}) \mgt\ldots\mgt
                             LM(\xi_k f_{j_k})$.
    \end{itemize}
\end{theorem}
\proof The standard proof applies to
{\it   i)}$\Rightarrow${\it  ii)},
{\it iii)}$\Rightarrow${\it  iv)},
{\it  iv)}$\Rightarrow${\it   v)} and
{\it   v)}$\Rightarrow${\it   i)}.

{\it  ii)}$\Rightarrow${\it iii)}:
Assume that {\it ii)} holds and let $f\in\Ktx$ homogeneous.
First we show the following\\
{\bf Claim:} If $f$ has a representation of the form
\begin{equation*}
    f = \sum_{i=1}^k\xi_i f_{j_i}\quad\text{with}\quad
    LM(f) = LM(\xi_1 f_{j_1}) \mgth LM(\xi_i f_{j_i}),\quad
    2\leq i\leq k,
\end{equation*}
then $\tf=f-\xi_1 f_{j_1}$ has a representation of the form
\begin{equation*}
    g\tf = \sum_{i=1}^{\tk}\txi_i f_{\tj_i}\quad\text{with}\quad
    LM(g\tf) = LM(\txi_1 f_{\tj_1})\mgth LM(\txi_i f_{\tj_i}),\quad
    2\leq i\leq \tk,
\end{equation*}
for a homogeneous $g\in \Ktx$ with $LM(g|_{t=1})=1$.

Indeed, after a permutation of summands $\tf=\sum_{i=2}^k\xi_i f_{j_i}$
with $LM(\tf)=LM(\xi_2 f_{j_2})=\ldots=LM(\xi_l f_{j_l})\mgth
LM(\xi_i f_{j_i})$, $l+1\leq i\leq k$ for a $l\in\{1,\ldots,k\}$.
If $l=1$ we are done. Otherwise let $c_i=LC(\xi_i f_{j_i})$, then
\begin{align*}
    \sum_{i=2}^l\xi_i f_{j_i} &=
            \underset{\text{lead terms cancel}}{\underbrace{\xi_2f_{j_2}-\frac{c_2}{c_3}\xi_3f_{j_3}}}+
            \left(1+\frac{c_2}{c_3}\right)\xi_3f_{j_3}+
            \sum_{i=4}^l\xi_i f_{j_i} \\
    &=      \eta h +\left(1+\frac{c_2}{c_3}\right)\xi_3 f_{j_3}+
            \sum_{i=4}^l\xi_i f_{j_i}
\end{align*}
for some $\eta\in\sT_t$ and $h\in\Spoly(f_{j_2},f_{j_3})$.
Since {\it ii)} holds $\eta h$ has a representation
$g\eta h=\sum_{i'=1}^{k'}\xi_i'f_{j_i'}$
with
$LM(g\xi_2f_{j_2})\mgth LM(g\eta h)=LM(\xi_1'f_{j_1'})\mgth
LM(\xi_i'f_{j_i'})$, $2\leq i\leq k'$. Then
\begin{align*}
    g\tf &= 
          \left(1+\frac{c_2}{c_3}\right)g\xi_3f_{j_3} +
          \sum_{i=4}^kg\xi_i f_{j_i}+
          \sum_{i=1}^{k'}\xi_i'f_{j_i'}
\end{align*}
and $LM(g\tf)=LM(g\xi_3f_{j_3})=\ldots=LM(g\xi_l f_{j_l})$ and
$LM(g\tf)\mgth LM(\xi_igf_{j_i})$, $l+1\leq i\leq k$ and
$LM(g\tf)\mgth LM(g\xi_i'f_{j'_i})$, $1\leq i\leq k'$. Thus (after expansion)
we have found a representation of $g\tf$ with $l-1$
lead terms instead of $l$. Then the claim follows by induction.

Now let $\bxa\in\sM$ and let $h\in\Spoly(f_i,f_j)$. Then there
exists a $g\in \Ktx$ with lead monomial 1 such that
$g\bxa h$ has a representation as in the claim.
Then the claim shows that during the computation of 
$\NormalForm(\bxa h,G)$ on never ends up
with an element which cannot be reduced except zero.
Thus  $\NormalForm(\bxa h,G)=0$ and {\it iii)} follows.
\proofend
\begin{corollary}
    Let $\mlt$ be such that
    \NormalForm\ terminates and reducing sets exist. If $G$ is
    a standard basis of the ideal $I\subseteq\Kx$, then 
    $G$ generates $\Iloc$ over $\Kxloc$.
\end{corollary}
%
%
\section{The Newton ordering}\label{sect:newton}
Let $\sL$ be a nonempty finite set of nonzero linear forms
$l\colon\RR^n\rightarrow\RR$ and
let $\bdelta\in(\RR_0^+)^n$ be fixed. We call the real number
\begin{align*}
    w_\bdelta(\balpha) &=
        \min\{l(\balpha+\bdelta)\mid l\in L\}
\end{align*}
the {\em weight} of $\balpha\in(\RR_0^+)^n$ with respect to $\sL$.
Let $C_\bdelta(l)=\{\balpha\in(\RR_0^+)^n\mid w_\bdelta(\balpha)=
l(\balpha+\bdelta)\}$, $l\in\sL$. A subscript of $\bdelta=\bzero$
will be omitted. Then $C_\bdelta(l)$ is a cone with respect to
$C(l)$: for all $\balpha\in C_\bdelta(l)$, $\bbeta\in C(l)$
by definition $w_\bdelta(\balpha+\bbeta)=l(\balpha+\bbeta)$,
so $\balpha+\bbeta\in C_\bdelta(l)$.
\begin{definition}
    $\sL$ is called {\em rational} if all cones $C(l)$, $l\in\sL$
    are rational (i.e.~every cone is the locus where a finite set of
    linear forms with rational coefficients is nonnegative).
\end{definition}
For the rest of this section let $\sL$ be rational and
$\bdelta\in(\QQ_0^+)^n$. Restricting
$w_\bdelta(.)$ to $\NN^n$ and identifying $\sM$ with $\NN^n$
we call
\begin{equation*}
    w_\bdelta(\bxa) = \min\{l(\balpha+\bdelta)\mid l\in\sL\}
\end{equation*}
the {\em Newton weight} of $\bxa\in\sM$. Define the Newton
weight of $f=\sum_{\balpha\in A}c_\balpha\bxa$ to be
\begin{equation*}
    w_\bdelta(f) =\min\{w_\bdelta(\bxa)\mid\balpha\in A\}.
\end{equation*}
The induced filtration
$\Kx_s=\{f\in \Kx\mid
w_\bdelta(f)\geq s\}$ of $\Kx$ is called a
{\em Newton filtration}. This filtration does not
distinguish all monomials in general. So take 
any semigroup order $\mlt_0$ and define
\begin{align*}
    \bxa\mlt\bxb\Longleftrightarrow &
     \ w_\bdelta(\bxa) > w_\bdelta(\bxb)\quad\text{or}\\
    &\ w_\bdelta(\bxa) = w_\bdelta(\bxb)\quad\text{and}\quad
       \bxa\mlt_0\bxb.
\end{align*}
for $\bxa$, $\bxb\in\sM$. This ordering is called a {\em Newton ordering}.
A Newton filtration (resp.~a Newton ordering) is a filtration
(resp.~an ordering) which arises in the above way.

\begin{figure}[H]
    \centering
    \setlength{\unitlength}{1mm}
    \begin{picture}(55,55)(0,0)
        \put(0,0){\epsfig{file=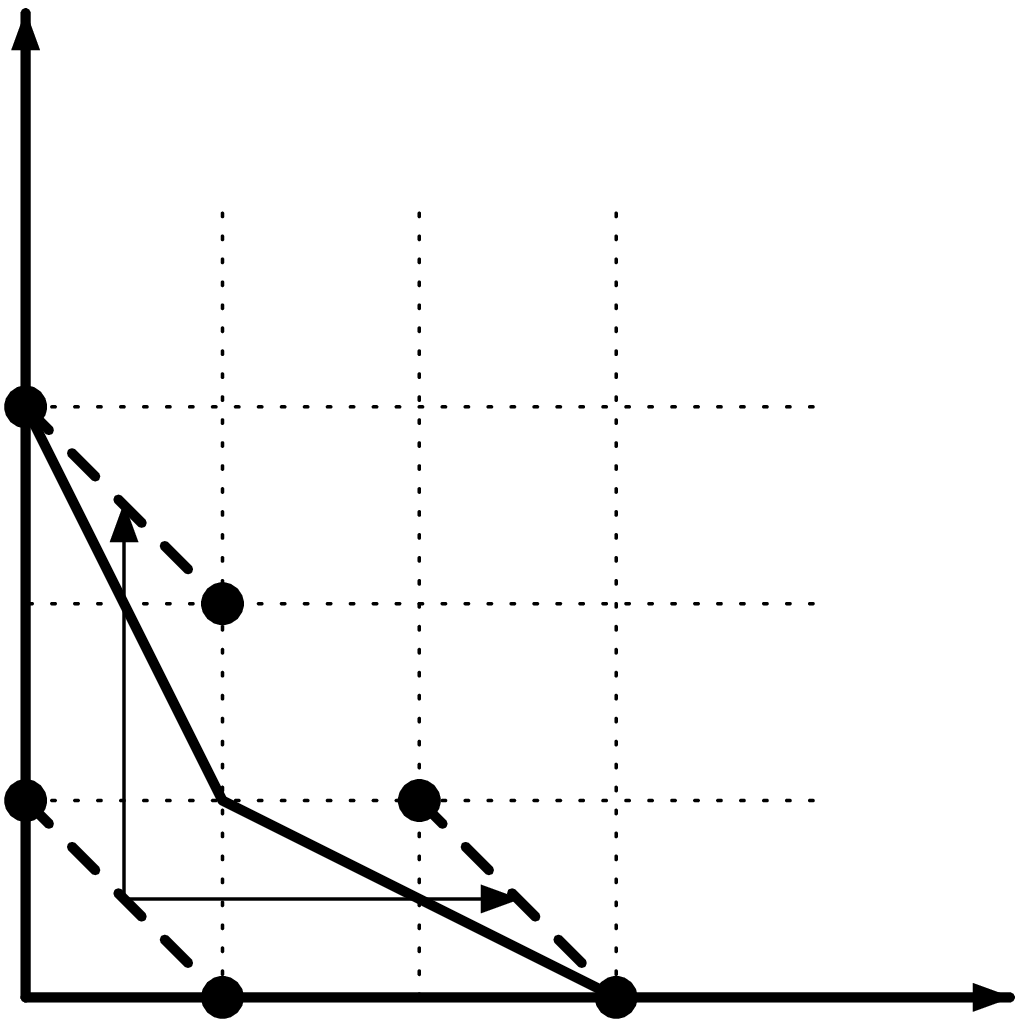,height=5cm}}
        \put(52,0){$x$}
        \put(0,52){$y$}
    \end{picture}
    \caption{$\sL=\{x+2y,2x+y\}$}
    \label{figure:newta}
\end{figure}

\begin{lemma}\label{lemma:newtonnormal}
    Let $\mlt$ be a Newton ordering.
    Then $\mlt$ is normal in the following three cases:
    \begin{itemize}
        \item[1)] All linear forms $l\in\sL$ have nonnegative coefficients.
        \item[2)] All linear forms $l\in\sL$ have nonpositive coefficients.
        \item[3)] $\bdelta=0$.
    \end{itemize}
\end{lemma}
\proof The first two cases are obvious.
So let $\bdelta=\bzero$ and $\bxa$, $\bxb\in\sM$.
We have to show that $\bxb\mlt 1$ implies $\bxa\bxb\mlt\bxa$.
There exist linear forms $l$, $l'$, $l''\in\sL$ such that
$w(\bxa)=l(\balpha)$,
$w(\bxb)=l'(\bbeta)$ and
$w(\bxa\bxb)=l''(\balpha+\bbeta)$. Therefore
\begin{align*}
    w(\bxa\bxb) &= l''(\balpha+\bbeta) =
                   l''(\balpha)+l''(\bbeta) \\
                &\geq l(\balpha)+l'(\bbeta) = w(\bxa)+w(\bxb).
\end{align*}
Then $\bxb\mlt 1$ implies either $w(\bxb)< 0\Rightarrow
w(\bxa\bxb)<w(\bxa)$ or
$w(\bxb)=0$ and $\bbeta\mlt_0 1\Rightarrow
w(\bxa\bxb)=w(\bxa)$ and $\bxa\bxb\mlt_0\bxa$. In both cases we get
$\bxb\bxa\mlt\bxa$.
\proofend
Note that $\mlt$ is an ordering on $\sM$ which is
not a semigroup ordering in general.
For example let $\sL=\{x+2y,2x+y\}$ and $f=x+y$ (figure \ref{figure:newta}).
Then $w_\bzero(x^2f)=w_\bzero(x^3+x^2y)=w(x^3)=3/2$, so
$LM(x^2f)=x^3$. On the other hand
$w_\bzero(y^2f)=w_\bzero(xy^2+y^3)=w(y^3)=3/2$, hence
$LM(y^2f)=y^3$. So $\mlt$ is not a semigroup order.
It follows from \cite{robbiano} that every semigroup ordering is
a Newton ordering with $|\sL|=1$ and $\bdelta=\bzero$.

If the coefficients of every $l\in L$ are positive (negative),
then $\mlt$ is a local (global) ordering.
Consider the set of monomials
\begin{equation*}
    \sM_l =\{\bxa\in\sM\mid w_\bdelta(\bxa) =l(\balpha+\bdelta)\},
    \quad l\in\sL.
\end{equation*}
Clearly $\sM_\bdelta(l)$ is a cone with respect to $\sM(l)$:
if $\bxa\in\sM_\bdelta(l)$ and $\bxa\in\sM(l)$,
then $\bxa\bxb\in\sM_\bdelta(l)$. Let
$R_\bdelta(l)$ be the $K$-algebra generated by the monomials
in $\sM_\bdelta(l)$. Then $R_\bdelta(l)$ is a 
$R(l)$-module. Since $\sL$ is rational, $R(l)$ is noetherian
and $R_\bdelta(l)$ is finitely generated since
$\bdelta$ has rational coordinates. Now let
$G=\fset\subset\Kx$ be a finite set of polynomials and
$I\subseteq\Kx$ an ideal. Define
$L(G)_\bdelta(l)=L(G)\cap\sM_\bdelta(l)$ and
$L(I)_\bdelta(l)=L(I)\cap\sM_\bdelta(l)$. Let 
$R(G)_\bdelta(l)$ (resp.~$R(I)_\bdelta(l)$) be the $K$-algebra
generated by the monomials in $L(G)_\bdelta(l)$ 
(resp.~$L(I)_\bdelta(l)$).
\begin{example}
    Let $\sL=\{l_1,l_2,l_3\}=\{3x+3y,2x+6y,6x+2y\}$
    and $\bdelta=(1,1)$. Then (figure \ref{figure:newtb})
    \begin{align*}
        R(l_1) &= \CC[xy,x^2y,x^3y,xy^2,xy^3],\quad  \\
        R(l_2) &= \CC[xy^3,y]\quad\text{and}\quad \\
        R(l_3) &= \CC[x,x^3y].
    \end{align*}
    \begin{figure}[H]
        \centering
        \setlength{\unitlength}{1mm}
        \begin{picture}(55,55)(0,0)
            \put(0,0){\epsfig{file=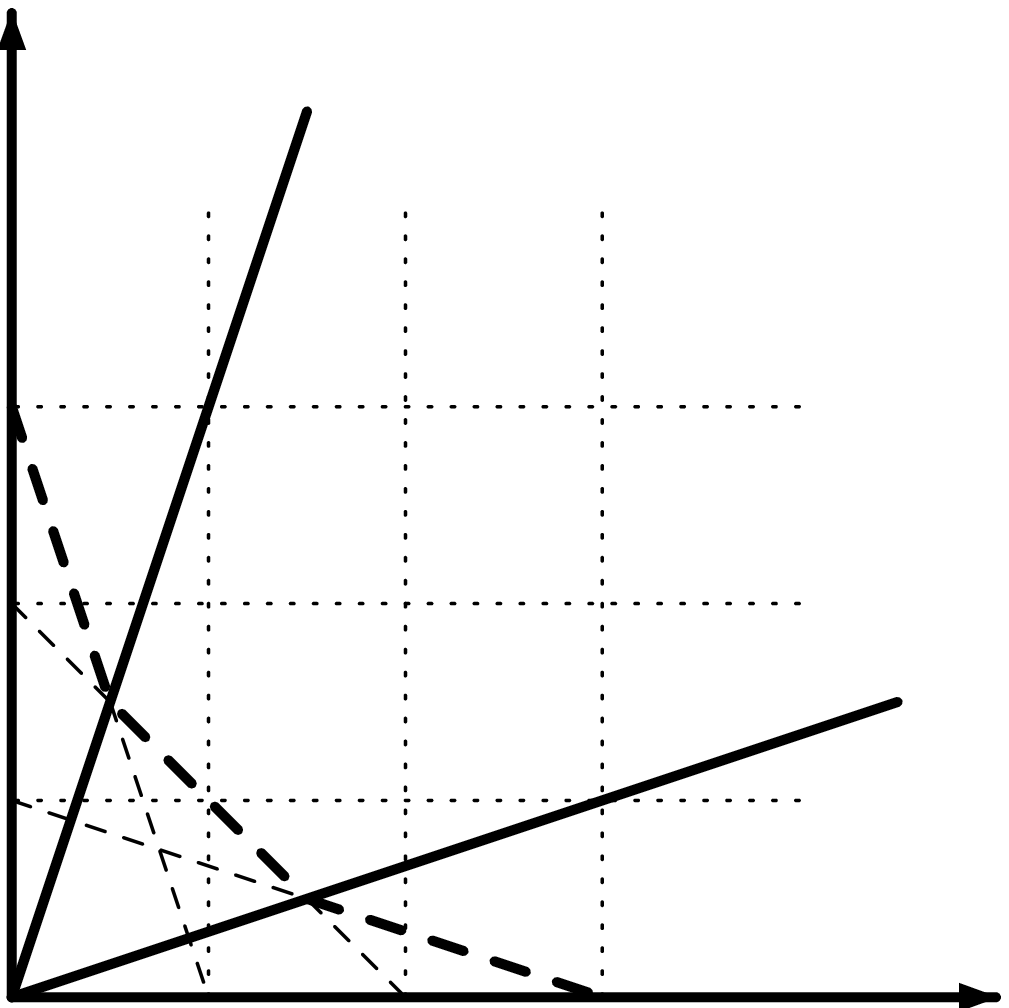,height=5cm}}
            \put(52,0){$x$}
            \put(0,52){$y$}
            \put(35,35){$\sM(l_1)$}
            \put(2,42){$\sM(l_2)$}
            \put(38,3){$\sM(l_3)$}
        \end{picture}
        \caption{$\sL=\{3x+3y,2x+6y,6x+2y\}$}
        \label{figure:newtb}
    \end{figure}
    The corresponding modules are (figure \ref{figure:newtc}):
    \begin{align*}
        R_\bdelta(l_1) &= x^2R(l_1)\oplus
                          x  R(l_1)\oplus
                             R(l_1)\oplus
                          y  R(l_1)\oplus
                          y^2R(l_1),\\
        R_\bdelta(l_2) &= y^2R(l_2)\quad\text{and}\quad\\
        R_\bdelta(l_3) &= x^2R(l_3).
    \end{align*}
    \begin{figure}[H]
        \centering
        \setlength{\unitlength}{1mm}
        \begin{picture}(55,55)(0,0)
            \put(0,0){\epsfig{file=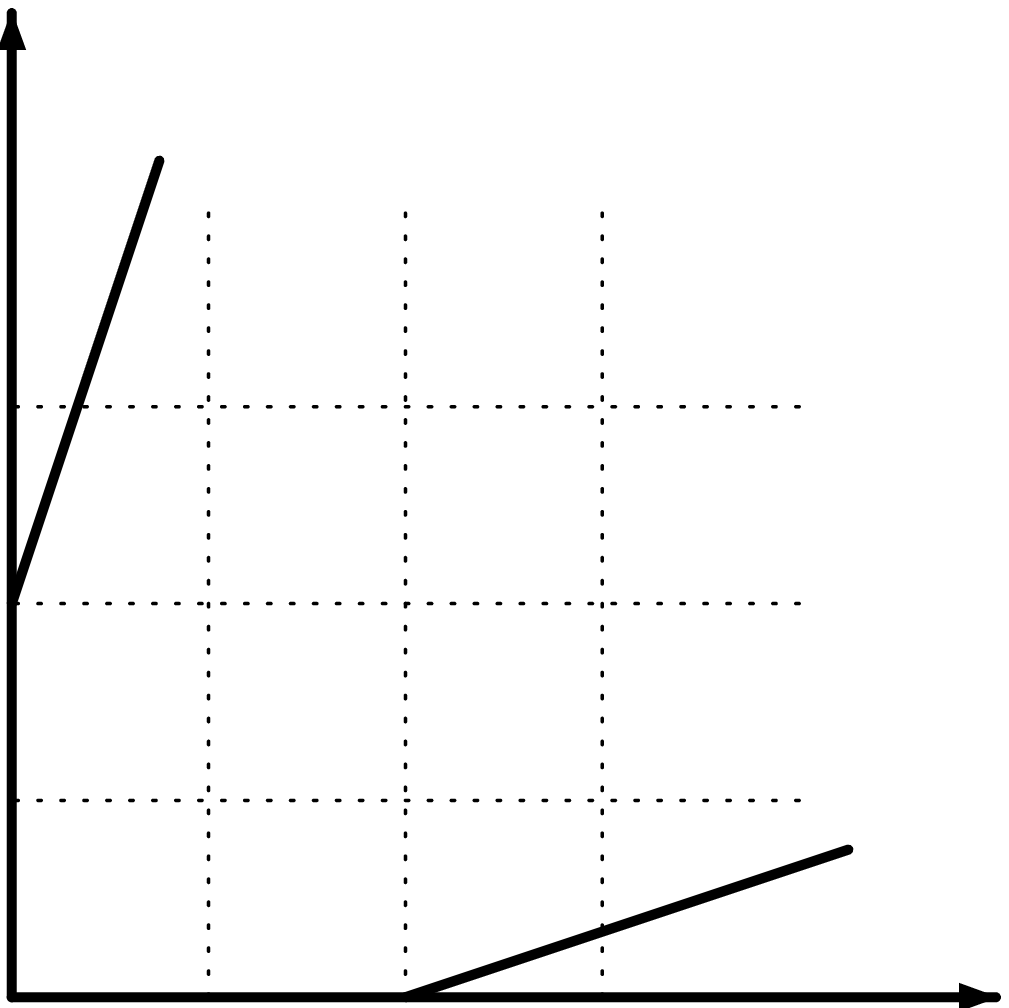,height=5cm}}
            \put(52,0){$x$}
            \put(0,52){$y$}
            \put(35,35){$\sM_\bdelta(l_1)$}
            \put(2,44){$\sM_\bdelta(l_2)$}
            \put(40,3){$\sM_\bdelta(l_3)$}
        \end{picture}
        \caption{$\sL=\{3x+3y,2x+6y,6x+2y\}$, $\bdelta=(1,1)$}
        \label{figure:newtc}
    \end{figure}
\end{example}
\begin{lemma}\label{lemma:sub}
    $L(G)_\bdelta(l)$ and $L(I)_\bdelta(l)$ are cones with respect
    to $\sM(l)$. $R(G)_\bdelta(l)$ and $R(I)_\bdelta(l)$ are
    finitely generated $R(l)$-submodules of $R_\bdelta(l)$.
\end{lemma}
\proof Let $f=\sum_{\balpha\in A} c_{\balpha}\bxa\in\Kx$ with
$LM(f)=\bx^{\balpha_0}\in\sM(I)_\bdelta(l)$. For every $\bxb\in\sM(l)$ and
every $\balpha\in A$ there exists a $l'\in\sL$ such that
\begin{align*}
    w_\bdelta(\bxa\bxb) &= l'(\balpha+\bbeta+\bdelta) =
                            l'(\balpha+\bdelta)+l'(\bbeta) \\
                            & \geq l(\balpha_0+\bdelta)+l(\bbeta) =
                            w_\bdelta(\bx^{\balpha_0}\bxb).
\end{align*}
so $LM(\bxb f)=\bxb LM(f)$.
Therefore both $L(G)_\bdelta(l)$ and $L(I)_\bdelta(l)$ are cones with respect
to $\sM(l)$. This shows that $R(G)_\bdelta(l)$ and $R(I)_\bdelta(l)$
are indeed $R(l)$-submodules of $R_\bdelta(l)$. They are finitely generated
since $R_\bdelta(l)$ is a noetherian module.
\proofend
\begin{proposition}\label{proposition:newtonnoether}
    Let $\mlt$ be a Newton ordering. Then $\mlt$ is noetherian.
\end{proposition}
\proof Let $L_0\subseteq L_1\subseteq L_2\subseteq\ldots$ be an
increasing sequence of lead monomial sets, i.e.~every $L_i$
is either $L(G)$ for a finite set of polynomials $G\subset\Kx$ or
$L(I)$ for an ideal $I\subseteq\Kx$. Let $l\in\sL$. By lemma \ref{lemma:sub}
the monomials of $L_i\cap\sM_\bdelta(l)$ span over $K$ a $R(l)$-submodule
$R_i$ of $R_\bdelta(l)$. Hence $R_0\subseteq R_1\subseteq R_2\subseteq\ldots$
is a increasing sequence of submodules. But $R_\bdelta(l)$ is
noetherian, so this sequence gets stationary and so does the
sequence $(L_i\cap\sM_\bdelta(l))_{i\in\NN}$. Since $|\sL|<\infty$
also the sequence $(L_i)_{i\in\NN}$ gets stationary.
\proofend
\begin{corollary}
    Let $\mlt$ be a Newton ordering.
    Then every ideal $I\subseteq K[\bx]$ has a standard basis
    with respect to $\mlt$.
\end{corollary}
\begin{corollary}
    Let $\mlt$ be a normal Newton ordering.
    Then \NormalForm\ terminates.
\end{corollary}
Let $G=\fset\subseteq\Ktx$ be a finite set of
homogeneous polynomials, let $I=(\fele)$ and $f\in I$ homogeneous.
We want to construct a reducing set for $(f,G)$.
Let $S_l$ be a finite set of monomials which generates
$R(\{f|_{t=1}\})_\bdelta(l)$ over $R(l)$, $l\in\sL$ and let
$S=\bigcup_{l\in\sL}S_l$. For every $\bxb\in S$ there exists
a monomial $\bxtb\in\sM$ such that $LM(\bxtb f|_{t=1})=\bxb$.
\begin{proposition}\label{proposition:newtonred}
    $R=\{\bxtb\mid\bxb\in S\}$ is a reducing set for $(f,G)$.
\end{proposition}
\proof Let $\bxa\in\sM$ and let
$t^{\gamma_0}\bxg=LM(\bxa f)\in\sM_\bdelta(l)$.
We have to find a representation of $\bxa f$ as in definition 
\ref{definition:reducing_set}. Then there exist monomials
$\bxtb\in R\cap\sM_\bdelta(l)$ and $\bxt\in\sM(l)$ such that
$\bxg=\bxb\bxt$ and $\bxa=\bxtb\bxt$. It follows from
lemma \ref{lemma:normal_form} that there exists a
homogeneous $g\in\Ktx$ with $LM(g|_{t=1})=1$,
$\xi_1,\ldots,\xi_k\in\sT_t$ and $h_1,\ldots,h_k\in G\cup
\{\NormalForm(\bxtb f,G)\mid \bxtb\in R\}$ such that
\begin{equation*}
    g\bxtb f=\sum_{i=1}^k \xi_i h_i,\qquad
    LM(g\bxtb f)=LM(\xi_1h_1)\mgt LM(\xi_ih_i),\quad 2\leq i\leq k.
\end{equation*}
But then $g\bxa f=\bxt g\bxtb f=\sum_{i=1}^k\bxt\xi_i h_i$ with
\begin{align*}
    LM(g\bxa f) &= LM(\bxt g\bxtb f) = \bxt LM(g\bxtb f)\\
                &= \bxt LM(\xi_1h_1)= LM(\bxt\xi_1 h_1)\\
                &\mgt LM(\bxt\xi_i h_i),
                \quad 2\leq i\leq k.
\end{align*}
This completes the proof.
\proofend
\begin{corollary}
    Let $\mlt$ be a normal Newton ordering.
    Then reducing sets do exist.
\end{corollary}
\begin{corollary}
    Let $\mlt$ be a normal Newton ordering.
    Then {\bf StandardBasis} terminates.
\end{corollary}
\section{Local orderings and zero dimensional ideals}\label{sect:local}
Let $\mlt$ be a local ordering and let $I\subseteq K[\bx]$ be a
zero dimensional ideal, i.e.~an ideal with $\dim_K(K[\bx]/I)<\infty$.
In this case $\Kxloc=\Kx_{(\bx)}=\Kx_{(x_1,\ldots,x_n)}$.
We show that there always exists a finite standard basis of $I$. Moreover
such a standard basis can be computed from a given standard basis
with respect to another local semigroup ordering using only linear algebra.
\begin{proposition}\label{proposition:always}
    Every zero dimensional ideal $I\subseteq K[\bx]$ has a
    standard basis with respect to any local order.
\end{proposition}
\proof Since $\dim I=0$ there
exists a $m\in\NN$ such that $(\bx)^m\subseteq I_{(\bx)}$.
Then $L(I)\setminus (\bx)^m=\{\bx^{\balpha_1},
\ldots,\bx^{\balpha_k}\}$
is a finite set of monomials. There exist polynomials
$f_1,\ldots,f_k\in\oI$ with $LM(f_i)=\bx^{\balpha_i}$,
$1\leq i\leq k$. Clearly $\{f_1,\ldots,f_k\}$ together
with all monomials of degree $m$ is a standard basis
of $I$ with respect to $\mlt$.
\proofend
For a local Newton ordering a standard basis can be computed
with our variant of Buchbergers algorithm. For an arbitrary local
ordering however it is not clear that this works. In particular, we
do not know the corresponding reducing sets.
One the other hand in the case of
zero dimensional ideals and global semigroup orderings, one can
do the following: given a standard basis with respect to a given
ordering, the FGLM-algorithm \cite{fglm} computes a standard basis
with respect to another ordering using only linear algebra.
The key is that zero dimensional ideals admit finite
defining systems.

There are two reasons why the FGML-algorithm will not work in our case:
Firstly, for a local ordering, a standard basis of $I$ does
not necessarily generate $\oI$ as a $K[\bx]$-module. 
The second reason is that the FGLM-algorithm depends on the fact that the set
of lead monomials has the structure of a monomial ideal, which is
in general only true for semigroup orderings.

Nevertheless it is possible to convert standard bases between
different local orderings.
Let $m\in\NN$ such that $(\bx)^m\subseteq
I_{(\bx)}$ and let $G=\fset$ be a standard basis of $I$ with respect
to $\mlt$.
Then $\oI$ is an ideal in
$K[\bx]$ with $(\bx)^m\subseteq I_{(\bx)}$. The following lemma
shows the relation between $G$ and $\oI$.
\begin{lemma}
    $G$ together with all monomials of degree $m$ generates
    $\oI$ as $K[\bx]$-module.
\end{lemma}
\proof This follows immediately from $L(I)=L(\oI)$.
\proofend
\subsection*{The algorithm}
First choose a monomial ideal $Z\subseteq\oI$ which contains a
power of the maximal ideal. The best choice for $Z$ would be
the maximal monomial ideal contained in $\oI$, but this ideal may
be expensive to compute. Since $Z$ contains
a power of the maximal ideal, $M\setminus Z$ is a finite set of
monomials. Consider the following partition of $M\setminus Z$:
\begin{align*}
    N &= \{\bx^\balpha\in M\setminus Z\mid \bx^\balpha\not\in L(I)\},\\
    H &= \{\bx^\balpha\in M\setminus Z\mid \bx^\balpha\in L(I)\}.
\end{align*}
Let $H=\{\bx^{\balpha_1},\ldots,
\bx^{\balpha_k}\}$ and $N=\{\bx^{\balpha_{k+1}},\ldots,
\bx^{\balpha_{k+\mu}}\}$. There exist polynomials
$h_1,\ldots,h_k\in\oI$ such that $LT(h_i)=\bx^{\balpha_i}$,
$1\leq i\leq k$. These polynomials are straightforward
to compute from $G$.
The ideal $\oI\subseteq\Kx$ is an infinite dimensional $K$-vector
space with basis $\{h_1,\ldots,h_k\}\cup(\sM\cap Z)$.
So after some linear algebra we can assume
$h_i-\bx^{\balpha_i}\in K\left<N\right>$ (this simplification
is just of cosmetic nature). Then
\begin{equation*}
    \left(\begin{array}{ccc|ccc}
    1 & & 0 & \hphantom{0} & \hphantom{0} & \hphantom{0}\\
      &\ddots & & \hphantom{0} & \ast & \hphantom{0}\\
    0 & & 1 & \hphantom{0} & \hphantom{0} & \hphantom{0}
    \end{array}\right)
    \begin{pmatrix}
    \bx^{\balpha_1}\\ 
    \vdots\\
    \bx^{\balpha_k}\\ 
    \bx^{\balpha_{k+1}}\\ 
    \vdots\\
    \bx^{\balpha_{k+\mu}}
    \end{pmatrix} =
    \begin{pmatrix}
    h_1\\
    \vdots\\
    h_k
    \end{pmatrix}.
\end{equation*}
Denote the $k\times(k+\mu)$ matrix on the left hand side by $A$.
Now reorder the monomials $\bx^{\alpha_1},\ldots,\bx^{\balpha_{k+\mu}}$
with respect to a second local ordering $\mlt_2$.
This defines a matrix $\tA$ (which arises from $A$ by permuting
columns) on which 
we perform Gauss elimination. Call the resulting matrix
again $\tA$. Then there exists a permutation
$\sigma\in S_{k+\mu}$ (operating on the columns of $\tA$) such that
\begin{equation*}
    \sigma(\tA)
    \begin{pmatrix}
    \bx^{\balpha_{\sigma(1)}}\\ 
    \vdots\\
    \bx^{\balpha_{\sigma(k)}}\\ 
    \bx^{\balpha_{\sigma(k+1)}}\\ 
    \vdots\\
    \bx^{\balpha_{\sigma(k+\mu)}}
    \end{pmatrix} =
    \left(\begin{array}{ccc|ccc}
    1 & & 0 & \hphantom{0} & \hphantom{0} & \hphantom{0}\\
      &\ddots & & \hphantom{0} & \ast & \hphantom{0}\\
    0 & & 1 & \hphantom{0} & \hphantom{0} & \hphantom{0}
    \end{array}\right)
    \begin{pmatrix}
    \bx^{\balpha_{\sigma(1)}}\\ 
    \vdots\\
    \bx^{\balpha_{\sigma(k)}}\\ 
    \bx^{\balpha_{\sigma(k+1)}}\\ 
    \vdots\\
    \bx^{\balpha_{\sigma(k+\mu)}}
    \end{pmatrix} =
    \begin{pmatrix}
    h'_1\\
    \vdots\\
    h'_k
    \end{pmatrix}
\end{equation*}
with $LT(h_i')=\bx^{\balpha_{\sigma(i)}}$ for $1\leq i\leq k$, thus
$\bx^{\balpha_{\sigma(1)}},\ldots\bx^{\balpha_{\sigma(k)}}\in L_{\mlt_2}(I)$.
As a $K$-vector space $\oI$ still has as basis
$\{h'_1,\ldots,h'_k\}\cup (Z\cap\sM)$, so
$\bx^{\balpha_{\sigma(k+ 1 )}},\ldots,
 \bx^{\balpha_{\sigma(k+\mu)}}\not\in L_{\mlt_2}(I)$.
Let $Z_{gen}$ be a finite set of generators of
the monomial ideal $Z$.
Then by construction 
$\{h'_1,\ldots,h'_k\}\cup Z_{gen}$
is a standard basis of $I$ with respect to $\mlt_2$.

Now we formulate the algorithm.
So let $G$ be a standard basis of $I$ with respect to $\mlt$.
Assume that we are given the following procedures:
\begin{description}
    \item[$\ComputeZgen(G,\mlt)$]\ \\
          Input: $G$ is a standard basis of a zero dimensional  ideal
          $I\subseteq\Kx$ with respect to the local ordering $\mlt$.\\
          Output: A reduced set of generators of a
          monomial ideal contained in $\oI$ which contains a power
          of $(\bx)$.
    \item [$\GaussEliminate(A)$]\ \\
          Input: $A$ is a $k\times(k+\mu)$ matrix over $K$ of
          rank $k$.\\
          Output: The $k\times(k+\mu)$ matrix over $K$ which
          arises from $A$ by Gaussian elimination.
\end{description}
\noindent Then the description of our algorithm is as follows.
\begin{algorithm}\label{algorithm:convert}\ \newline\noindent
    $G_2:=\StandardBasisChange(G,\mlt,\mlt_2)$\\
    \hspace*{1cm}$Z_{gen}:=\ComputeZgen(G,\mlt)$\\
    \hspace*{1cm}$M:=\{\bx^\balpha\in\sM\mid\bxa\not\in (Z_{gen})\}$\\
    \hspace*{1cm}$N:=\{\bx^\balpha\in M\mid \bxa\not\in L(G)\}$\\
    \hspace*{1cm}$H:=M\setminus N$\\
    \hspace*{1cm}$F:=\emptyset$\\
    \hspace*{1cm}FOR ALL $\bxa\in H$ DO\\
    \hspace*{2cm}choose the first $f\in G$ with $LM(\eta f)=\bxa$ for
                 some $\eta\in\sT$\\
    \hspace*{2cm}$F:=F\cup\{\eta f\}$\\
    \hspace*{1cm}$A:=$ the matrix of coefficients of all $f\in F$ with
                 respect to $M$\\
    \hspace*{2cm}sorted according to $\mlt_2$\\
    \hspace*{1cm}$A:={\bf GaussEliminate}(A)$\\
    \hspace*{1cm}$F:=$ the elements of the vector which arises from
                 multiplying $A$ with\\ 
    \hspace*{2cm}the column vector of the elements of $M$ sorted according
                 to $\mlt_2$\\
    \hspace*{1cm}$G_2:=F\cup Z_{gen}$
\end{algorithm}
Altogether we have proved the
\begin{proposition}\label{proposition:convert}
    Let $G$ be a standard basis of the ideal $I\subseteq\Kx$
    with respect to the local ordering $\mlt$ and let $\mlt_2$
    be another local ordering. Then
    \StandardBasisChange\ computes a standard basis of $I$
    with respect to $\mlt_2$.
\end{proposition}
\begin{corollary}
    For every ideal $I\subseteq\Kx$ and every local ordering $\mlt$
    \begin{equation*}
        \dim_K\Kx/\oI = \left|\{\bxa\in\sM\mid\bxa\not\in L(I)\}\right|.
    \end{equation*}
\end{corollary}
For degree orderings like $degrevlex^-$ the set of generators $Z_{gen}$
can be computed as follows: One computes the highest corner
$\bxt$ of $I$, i.e.~the smallest monomial not contained
in $L(I)$. Then $\bxa\mlt\bxt$ implies $\bxa\in\oI$. For degree
orderings the set $\{\bxa\in\sM\mid\bxa\mget\bxt\}$ is
finite. Let
$Z_1=\{\bxa\in\sM\mid |\balpha|=|\btau|,\bxa\mlt\bxt\}$ and
$Z_2=\{\bxa\in\sM\mid |\balpha|=|\btau|+1,\bxa\not\in(Z_1)\}$.
Then we may take $Z_{gen}=Z_1\cup Z_2$.
This will also work (with a few modifications) in the case
of a weighted degree ordering or a Newton ordering.
\section{Applications}\label{sect:applications}
Let $\mlt$ be a normal Newton ordering, $I\subseteq\Kx$ an ideal
and $G$ a standard basis of $I$ with respect to $\mlt$. The Newton
filtration on $\Kx$ induces a filtration on the quotient 
$\Kx/\oI$. Let $Gr_\loc\Kx/\oI=\bigoplus_{s\in\RR}
(\Kx/\oI)_s/(\Kx/\oI)_{>s}$ be the associated graded
module. Now assume that $\dim_K(\Kx/\oI)_s/(\Kx/\oI)_{>s}<\infty$ for
all $s\in\RR$ and let $P(t)$ be the Poincare series of
$Gr_\loc\Kx/\oI$.
The coefficients of this series can be computed as follows.
\begin{proposition}\label{proposition:filtration}
    Let $s\in\RR$. The coefficient $c_s$ of $t^s$ in
    $P(t)$ is
    \begin{equation*}
        c_s = \dim_K(\Kx/\oI)_s/(\Kx/\oI)_{>s} = 
             \left|\left\{\bxa\in M\mid
               \bxa\not\in L(I), w_\bdelta(\bxa)=s
              \right\}\right|.
    \end{equation*}
\end{proposition}
\proof By definition the elements of
$B=\left\{\bxa\in\sM\mid\bxa\not\in L(I),
w_\bdelta(\bxa)=s\right\}$ are linearly
independent in the $K$-vector space
$V=\left(\Kx\right/\oI)_s/\left(\Kx\right/\oI)_{>s}$.
Now let $\bxa\in L(I)$ be the smallest monomial of Newton
weight $w_\bdelta(\bxa)=s$. Then there exists
a $f\in\oI$ with $LM(f)=\bxa$. By construction
$f-LT(f)\in K\left<B\right>\cup\Kx_{>s}$, so
$\bxa\in K\left<B\right>$ in $V$.
By induction this holds for
all $\bxa\in L(I)$ of Newton weight $s$.
\proofend
\subsection{The spectrum of a hypersurface singularity}
    %
    Let $f\colon (\CC^n,\bzero)\rightarrow(\CC,0)$ be the germ
    of a holomorphic function with an isolated critical point at
    $\bzero$. Let $J_f=
    \left(\partial f/\partial x_1,\ldots,
    \partial f/\partial x_n\right)\subseteq\CC\{\bx\}$ be the Jacobi
    ideal of $f$. The Milnor number of $f$
    \begin{equation*}
        \mu(f) =\dim_{\CC}\left(\CC\{\bx\}/J_f\right)
    \end{equation*}
    is a topological invariant of $f$. Since $f$ has an isolated 
    singularity $\mu(f)$ is finite. Then $f$ is
    $(\mu+1)$-determined, which means that $f$ does not change
    its analytical type if we forget about terms of order
    higher than $\mu+1$. So we can assume that $f$ is given by
    a polynomial. Moreover there exists a $m\in\NN_0$ such
    that $(\bx)^m\subseteq J_f$, so $\CC\{\bx\}/J_f\simeq\CC[\bx]_{(\bx)}/J_f$.
    This means that we can compute the Milnor number of $f$    
    from a standard basis of $J_f$ with respect to a local ordering
    $\mlt$ as
    \begin{equation*}
        \mu(f) =\left|\sM\setminus L(J_f)\right|.
    \end{equation*}
    The spectrum of $f$ is an analytical invariant of $f$ which is
    finer than the Milnor number. It consists of $\mu=\mu(f)$
    rational numbers $\{\alpha_1,\ldots,\alpha_\mu\}$ which
    often are written in the form $\sum_{i=1}^\mu t^{\alpha_i}$.
    These rational numbers are defined in terms of the
    mixed Hodge structure on the vanishing cohomology
    of the Milnor fibre of $f$ \cite{arnold}. In many cases the spectrum
    is determined by the Newton polygon of $f$ which is constructed
    as follows. If $f=\sum_{\balpha\in A}c_\balpha\bxa$ let
    $Or(\bx^\balpha)$ be the positive orthant in $\RR^n$ centred
    at $\balpha\in\NN_0^n\subset\RR^n$, $\balpha\in A$. Then
    the Newton polygon $\Sigma_f$ of $f$ is the convex hull of the union
    $\bigcup_{\balpha\in A}Or(\bx^\balpha)$.

    Let $\Sigma_f^c$ be the set of compact faces of $\Sigma_f$.
    For every compact face $\sigma\in\Sigma_f^c$, let
    $f_\sigma =\sum_{\balpha\in A\cap\sigma}c_\balpha\bxa$
    be the sum of the monomials of $f$ with support on $\sigma$.
    \begin{definition}
        $f$ has {\em nondegenerate principal part} if for every
        $\sigma\in\Sigma_f^c$ the polynomials
        $\partial f_\sigma/\partial x_1,\ldots,
        \partial f_\sigma/\partial x_n$ do not have a common
        root in $(\CC^\ast)^n$.
    \end{definition}
    For every compact face $\sigma\in\Sigma_f^c$ there exists a
    unique linear form $l_\sigma\colon\QQ^n\rightarrow\QQ$ with
    $\left.l_\sigma\right|_\sigma\equiv 1$.
    Let $\sL=\{l_\sigma\mid\sigma\in\Sigma_f^c\}$
    and let $\bdelta=(1,\ldots,1)$. Consider the Newton
    filtration on $\CC[\bx]$ associated to $\sL$ and $\bdelta$.
    This induces a filtration on $\CC\{\bx\}$. The following
    theorem is due to M.~Saito.
    \begin{theorem}\label{theorem:saito}
        \cite{saito,kv} If $f$ has nondegenerate principal part then the
        spectrum of $f$ coincides with the Poincare series
        of the artinian module $\CC\{\bx\}/J_f$ graded
        by the Newton filtration.
    \end{theorem}
    Now let $\mlt$ be a Newton ordering for $\sL$ and $\bdelta$.
    Clearly $\mlt$ is local. Since we have assumed that $f$
    is given by a polynomial, as a consequence of
    proposition \ref{proposition:filtration} and theorem \ref{theorem:saito}
    we get the
    \begin{corollary}\label{corollary:spectrum}
        If $f$ has nondegenerate principal part then
        the spectrum of $f$ consists of the $\mu$ rational
        numbers
        \begin{equation*}
            \{w_\bdelta\left(\bxa\right)\mid
            \bxa\not\in L(J_f)\}.
        \end{equation*}
    \end{corollary}
    \begin{example}
        Consider the surface singularity
        \begin{align*}
            f=&x^{12}+y^{12}+z^{12}+x^5y^5+x^5z^5+y^5z^5\\
             &+xyz(x^2y^2+x^2z^2+y^2z^2)+x^2y^2z^2.
        \end{align*}
        The principal part of $f$ is nondegenerate and the Newton
        polygon of $f$ has twelve faces:
        \begin{figure}[H]
            \centering
            \setlength{\unitlength}{1mm}
            \begin{picture}(100,100)(0,0)
            \put(10,100){\epsfig{file=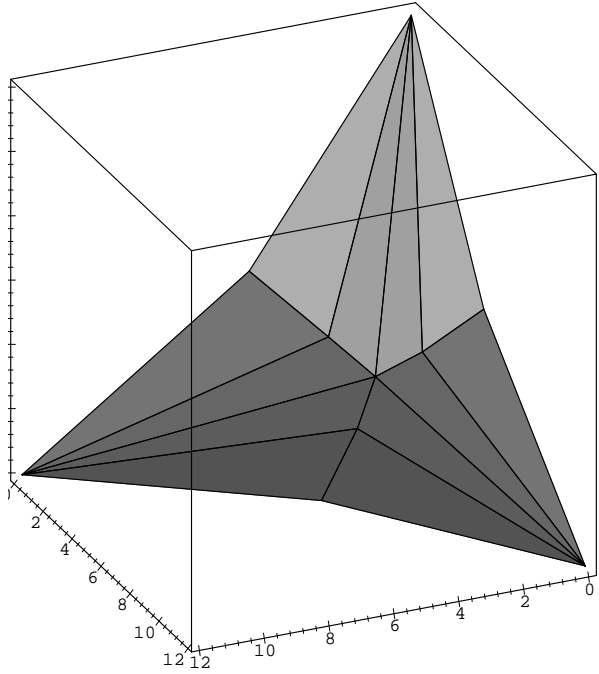,width=10cm,angle=270,%
                                clip=,
                                bbllx=50,
                                bblly=200,
                                bburx=550,
                                bbury=700}}
            \end{picture}
            \caption{Newton polygon of $f$}
            \label{figure:npoly}
        \end{figure}
        Using corollary \ref{corollary:spectrum} and algorithm
        \ref{algorithm:convert}, the spectrum of $f$ is
        \begin{align*}
            t^\frac{1}{2}&+
            3t^\frac{7}{12}+
            3t^\frac{2}{3}+
            6t^\frac{3}{4}+
            9t^\frac{5}{6}+
            9t^\frac{11}{12}+
            13t\\
            +&
            18t^\frac{13}{12}+
            18t^\frac{7}{6}+
            21t^\frac{5}{4}+
            24t^\frac{4}{3}+
            24t^\frac{17}{12}+
            25t^\frac{3}{2}\\
            +&
            24t^\frac{19}{12}+
            24t^\frac{5}{3}+
            21t^\frac{7}{4}+
            18t^\frac{11}{6}+
            18t^\frac{23}{12}+
            13t^2\\
            +&
            9t^\frac{25}{12}+
            9t^\frac{13}{6}+
            6t^\frac{9}{4}+
            3t^\frac{7}{3}+
            3t^\frac{29}{12}+
            t^\frac{5}{2}.
        \end{align*}
        In particular $\mu(f)=323$ and $p_g(f)=44$.
    \end{example}
    %
%
\section{Concluding remarks}
This work arose from the desire to have a program which computes
the spectrum of an isolated hypersurface singularity. Our variant of the
Buchberger algorithm has been successfully coded. However it turned out
that it is much faster to compute a standard basis with respect
to a semigroup ordering and to convert it to a standard basis
with respect to the corresponding Newton ordering using
$\StandardBasisChange$. This not really a surprise:
standard basis algorithms for semigroup orderings have been
carefully optimised for years, see \cite{ggmnppss,macaulay}.

The computation of the spectrum of an isolated hypersurface
singularity with nondegenerate principal part is now
implemented in the computer
algebra program {\tt Singular} and will hopefully be
available in further releases. For convenience, our web site
\begin{center}
    {\small\tt
    www.mathematik.uni-mainz.de/AlgebraischeGeometrie/Spectrum/index.shtml}
\end{center}
offers an interface to this implementation.
\subsection*{Acknowledgements}
I would like to thank D.~van Straten for telling me about
this problem, G.-M.~Greuel for valuable suggestions and
H.~Sch\"onemann for his support during the implementation.
G.-M.~Greuel, G.~Pfister and H.~Sch\"onemann gave me
the opportunity to implement the computation of the spectrum
of an isolated hypersurface singularity into {\tt Singular}.
\noindent%
Johannes Gutenberg-Universit\"at\\
Fachbereich 17\\
Staudinger-Weg 9\\
55099 Mainz, Germany\\
{\tt endrass@mathematik.uni-mainz.de}
\end{document}